\documentclass[11pt,reqno]{amsart}
\usepackage{amsmath,amssymb}

 \makeatletter
 \evensidemargin  \oddsidemargin
 \makeatother

\begin{document}
\thispagestyle{empty}
 \setcounter{page}{1}

\begin{center}
{\large\bf  Quadratic--quartic functional equations in RN--spaces

\vskip.20in

{\bf M. Bavand Savadkouhi } \\[2mm]

{\footnotesize Department of Mathematics,
Semnan University,\\ P. O. Box 35195-363, Semnan, Iran\\
[-1mm] e-mail: {\tt bavand.m@gmail.com}}

{\bf M. Eshaghi Gordji } \\[2mm]

{\footnotesize Department of Mathematics,
Semnan University,\\ P. O. Box 35195-363, Semnan, Iran\\
[-1mm] e-mail: {\tt madjid.eshaghi@gmail.com}}

{\bf  Choonkil Park} \\[2mm]

{\footnotesize Department of Mathematics,
Hanyang University,\\ Seoul 133-791, South Korea\\
[-1mm] e-mail: {\tt baak@hanyang.ac.kr}}}
\end{center}
\vskip 5mm

\noindent{\footnotesize{\bf Abstract.} In this paper, we obtain
the general solution and the stability result for the following
functional equation in random normed spaces (in the sense of
Sherstnev) under arbitrary $t$-norms
$$f(2x+y)+f(2x-y)=4[f(x+y)+f(x-y)]+2[f(2x)-4f(x)]-6f(y).$$
\vskip.10in 

  \newtheorem{df}{Definition}[section]
  \newtheorem{rk}[df]{Remark}
   \newtheorem{lem}[df]{Lemma}
   \newtheorem{thm}[df]{Theorem}
   \newtheorem{pro}[df]{Proposition}
   \newtheorem{cor}[df]{Corollary}
   \newtheorem{ex}[df]{Example}

 \setcounter{section}{0}
 \numberwithin{equation}{section}
\vskip .2in
\begin{center}
\section{Introduction}
\end{center}
The stability problem of functional equations originated from a
question of Ulam \cite{Ul} in $1940,$ concerning the stability of
group homomorphisms. Let $(G_1,.)$ be a group and let $(G_2,*,d)$ be
a metric group with the metric $d(.,.).$ Given $\epsilon >0$, does
there exist a $\delta>0$ such that if a mapping $h:G_1 \to G_2$
satisfies the inequality $d(h(x.y),h(x)*h(y)) <\delta$ for all
$x,y\in G_1$, then there exists a homomorphism $H:G_1 \to G_2$ with
$d(h(x),H(x))<\epsilon$ for all $x \in G_1?$ In the other words,
under what condition does there exists a homomorphism near an
approximate homomorphism? The concept of stability for functional
equation arises when we replace the functional equation by an
inequality which acts as a perturbation of the equation.  Hyers \cite{Hy} gave a first affirmative answer to the
question of Ulam for Banach spaces. Let $f:E \to E'$ be a mapping
between Banach spaces such that
$$\|f(x+y)-f(x)-f(y)\|\leq \delta $$
for all $x,y\in E$ and some $\delta>0.$ Then there exists a
unique additive mapping $T:E \to E'$ such that
$$\|f(x)-T(x)\|\leq \delta$$
for all $x\in E.$ Moreover, if $f(tx)$ is continuous in $t\in \Bbb
R$ for each fixed $x\in E,$ then $T$ is $\mathbb{R}$-linear. In $1978,$ Th. M.
Rassias \cite{TRa} provided a generalization of the Hyers' theorem
which allows the Cauchy difference to be unbounded. In $1991,$ Z.
Gajda \cite{Ga} answered the question for the case $p>1$, which
was raised by Rassias.  This new concept is known as
Hyers-Ulam-Rassias stability of functional equations (see
\cite{Ac-Dh,Ao,Bo,G,Hy-Is-Ra,Is-Ra,TRa1,TRa2}).
The functional equation
$$f(x+y)+f(x-y)=2f(x)+2f(y).\eqno(1.1)$$
is related to a symmetric bi-additive mapping. It is natural that
this equation is called a quadratic functional equation. In
particular, every solution of the quadratic functional  equation $(1.1)$ is
said to be a quadratic mapping. It is well known that a mapping
$f$ between real vector spaces is quadratic if and only if there
exits a unique symmetric bi-additive mapping $B$ such that
$f(x)=B(x,x)$ for all $x$ (see \cite{Ac-Dh,Kl}). The bi-additive
mapping $B$ is given by
$$B(x,y)=\frac{1}{4}(f(x+y)-f(x-y)).\eqno(1.2)$$
The Hyers-Ulam-Rassias stability problem for the quadratic functional
equation $(1.1)$ was proved by Skof for mappings $f:A \to B$, where
$A$ is a normed space and $B$ is a Banach space (see \cite{Sk}). Cholewa
\cite{Ch} noticed that the theorem of Skof is still true if relevant
domain $A$ is replaced an abelian group. In \cite{Cz},
Czerwik proved the Hyers-Ulam-Rassias stability of the functional equation
$(1.1)$. Grabiec \cite{Gr} has generalized these results
mentioned above.\\
In \cite{Pa-Ba}, W. Park and J. Bae considered the
following quartic functional equation
$$f(x+2y)+f(x-2y)=4[f(x+y)+f(x-y)+6f(y)]-6f(x).\eqno(1.3)$$
In fact, they proved that a mapping $f$ between two real vector
spaces $X$ and $Y$ is a solution of $(1.3)$ if and only if there
exists a unique symmetric multi-additive mapping $M: X^4 \rightarrow Y$ such that $f(x)=M(x,x,x,x)$ for all $x$.
It is easy to show that the function $f(x)=x^4$ satisfies the
functional equation $(1.3),$ which is called a quartic functional
equation (see also \cite{Ch-Sa}). In addition,  Kim \cite{Ki} has
obtained the  Hyers-Ulam-Rassias stability
for a mixed type of quartic and quadratic functional equation.\\
The Hyers-Ulam-Rassias stability of different
functional equations in random normed and fuzzy normed spaces has
been recently studied in \cite{Mi1}-\cite{Mi-Mo2}. It should be
noticed that in all these papers the triangle inequality is
expressed by using the strongest triangular norm $T_M$.

The aim of
this paper is to investigate the stability of the
additive-quadratic functional equation in random normed spaces (in
the sense of Sherstnev) under arbitrary
continuous $t$-norms.\\
In the sequel, we adopt the usual terminology, notations and
conventions of the theory of random normed spaces, as in
\cite{Ch-Ch-Ka,Mi2,Mi-Ra,Sc-Sk,Sh}. Throughout this paper,
$\Delta^{+}$ is the space of distribution functions that is, the
space of all mappings $F:\Bbb R \cup \{-\infty,\infty\} \to
[0,1]$ such that $F$ is left-continuous and non-decreasing on
$\Bbb R,$ $F(0)=0$ and $F(+\infty)=1.$ $D^{+}$ is a subset of
$\Delta^{+}$ consisting of all functions $F \in \Delta^{+}$ for
which $l^{-}F(+\infty)=1$, where $l^-f(x)$ denotes the left limit
of the function $f$ at the point $x$, that is, $l^{-}f(x)=\lim_{t
\to x^{-}}f(t)$. The space $\Delta^{+}$ is partially ordered by
the usual point-wise ordering of functions, i.e., $F \leq G$ if
and only if $F(t) \leq G(t)$ for all $t$ in $\Bbb R$. The maximal
element for $\Delta^{+}$ in this order is the distribution
function $\varepsilon_{0}$  given by
$$\varepsilon_{0}(t)=\begin{cases}
\begin{array}{ccc}
0, & \text{if} \; t\leq 0, \\
1, & \text{if} \; t>0 .\\
\end{array}
\end{cases}$$
\begin{df}(\cite{Sc-Sk}).
A mapping $T:[0,1] \times [0,1]\to [0,1]$ is a continuous
triangular norm (briefly, a continuous $t$-norm) if $T$ satisfies
the following conditions:\\
$(a)$ $T$ is commutative and associative;\\
$(b)$ $T$ is continuous;\\
$(c)$ $T(a,1)=a$ for all $a\in [0,1]$;\\
$(d)$ $T(a,b)\leq T(c,d)$ whenever $a\leq c$ and $b\leq d$ for all
$a,b,c,d\in [0,1]$.
\end{df}
Typical examples of continuous $t$-norms are $T_P(a,b)=ab$, $T_M
(a,b)=\min(a,b)$ and $T_{L}(a,b)=\max(a+b-1,0)$ (the Lukasiewicz
$t$-norm). Recall (see \cite{Ha-Pa, Ha-Pa-Bu}) that if $T$
is a $t$-norm and $\{x_n\}$ is a given sequence of numbers in
$[0,1]$, then $T_{i=1}^n x_i$ is defined recurrently by
$T^1_{i=1}x_i=x_1$ and $T^n_{i=1}x_i=T(T^{n-1}_{i=1}x_i)$ for
$n \geq 2.$ $T_{i=n}^\infty x_i$ is defined as $T_{i=1}^\infty
x_{n+i}.$ It is known (\cite{Ha-Pa-Bu}) that for the Lukasiewicz
$t$-norm the following implication holds:
$$ \lim_{n \to \infty}{(T_L)}_{i=1}^{\infty} x_{n+i}=1 \Longleftrightarrow
\sum_{n=1}^{\infty}(1-x_{n})<\infty.$$
\begin{df}(\cite{Sh}).
A \emph{random normed space} (briefly, RN-space) is a triple
$(X,\mu,T)$, where $X$ is a vector space, $T$ is a continuous
$t$-norm and $\mu$ is a mapping from $X$ into $D^{+}$ such that
the following conditions hold:\\
$(RN1)$ $\mu_x(t)=\varepsilon_{0}(t)$ for all $t>0$ if and only if $x=0$;\\
$(RN2)$ $\mu_{\alpha x}(t)=\mu_{x}(\frac{t}{|\alpha|})$ for all $x\in X$, $\alpha \neq 0$;\\
$(RN3)$ $\mu_{x+y}(t+s)\geq T(\mu_{x}(t),\mu_{y}(s))$ for all
$x,y\in X$ and $t,s \geq 0.$
\end{df}
Every normed space $(X,\|.\|)$ defines a random normed space
$(X,\mu,T_M)$, where $$\mu_x(t)=\frac{t}{t+\|x\|}$$ for all $t>0,$
and $T_M$ is the minimum $t$-norm. This space is called the
induced random normed space.
\begin{df}
Let $(X,\mu,T)$ be an RN-space.\\
$(1)$ A sequence $\{x_{n}\}$ in $X$ is  said to be
\emph{convergent} to $x$ in $X$ if, for every $\epsilon>0$ and
$\lambda>0$, there exists a positive integer $N$ such that
$\mu_{x_{n}-x}(\epsilon )>1-\lambda$ whenever $n\geq N$.\\
$(2)$ A sequence $\{x_{n}\}$ in $X$ is called a \emph{Cauchy
sequence} if, for every $\epsilon>0$ and $\lambda>0$, there exists
a positive integer $N$ such that
$\mu_{x_{n}-x_{m}}(\epsilon)>1-\lambda$
whenever $n \geq m \geq N$.\\
$(3)$ An RN-space $(X,\mu,T)$ is said to be \emph{complete} if and
only if every Cauchy sequence in $X$ is convergent to a point in
$X$.
\end{df}
\begin{thm}(\cite{Sc-Sk}).
If $(X,\mu,T)$ is an RN-space and $\{x_{n}\}$ is a sequence such
that $x_n \to x$, then $\lim_{n\to\infty} \mu_{x_n}(t)=\mu_{x}(t)$
almost everywhere.
\end{thm}
Recently, M. Eshaghi Gordji et al. establish the stability of cubic,
quadratic and additive-quadratic functional equations in RN-spaces
 (see \cite{Es-Ra-Ba} and \cite{E-R-B}).\\
 In this paper, we deal with the following functional equation
$$f(2x+y)+f(2x-y)=4[f(x+y)+f(x-y)]+2[f(2x)-4f(x)]-6f(y)\eqno(1.4)$$ on
RN-spaces. It is easy to see that the function
$f(x)=ax^4+bx^2$ is a solution of $(1.4).$\\

In Section 2, we investigate the general solution of the functional
equation $(1.4)$ when $f$ is a mapping between vector spaces and in
Section 3, we establish the stability of the functional equation
$(1.4)$ in RN-spaces. \vskip5mm
%================================================================
\section{General solution}
\setcounter{equation}{0}
%================================================================
We need the following lemma for solution of $(1.4).$ Throughout this
section $X$ and $Y$ are vector spaces.
\begin{lem}\label{t2}
If a mapping $f:X\longrightarrow Y$ satisfies $(1.4)$ for all $x,y
\in X,$ then $f$ is quadratic-quartic.
\end{lem}
\begin{proof}
We show that the mappings $g:X \longrightarrow Y$ defined by
$g(x):=f(2x)-16f(x)$ and $h:X \longrightarrow Y$ defined by
$h(x):=f(2x)-4f(x)$ are quadratic and quartic, respectively.

Letting $x=y=0$ in $(1.4),$ we have
$f(0)=0$. Putting $x=0$ in $(1.4)$, we get $f(-y)=f(y)$. Thus
the mapping $f$ is even. Replacing $y$ by $2y$ in $(1.4),$ we get
$$f(2x+2y)+f(2x-2y)=4[f(x+2y)+f(x-2y)]+2[f(2x)-4f(x)]-6f(2y) \eqno(2.1)$$
for all $x,y \in X$. Interchanging $x$ with $y$ in $(1.4),$ we
obtain
$$f(2y+x)+f(2y-x)=4[f(y+x)+f(y-x)]+2[f(2y)-4f(y)]-6f(x) \eqno(2.2)$$
for all $x,y \in X$. Since $f$ is even, by $(2.2),$ one gets
$$f(x+2y)+f(x-2y)=4[f(x+y)+f(x-y)]+2[f(2y)-4f(y)]-6f(x) \eqno(2.3)$$
for all $x,y \in X.$ It follows from $(2.1)$ and $(2.3)$ that
$$[f(2(x+y))-16f(x+y)]+[f(2(x-y))-16f(x-y)]=2[f(2x)-16f(x)]+2[f(2y)-16f(y)]$$
for all $x,y \in X$. This means that
$$g(x+y)+g(x-y)=2g(x)+2g(y)$$
for all $x,y \in X$. Therefore, the mapping $g:X \rightarrow Y$
is quadratic.

To prove that $h:X \rightarrow Y$ is quartic, we have to show that
$$h(x+2y)+h(x-2y)=4[h(x+y)+h(x-y)+6h(y)]-6h(x)$$
for all $x,y \in X.$ Since $f$ is even, the mapping $h$ is even.
Now if we interchange $x$ with $y$ in the last equation, we get
$$h(2x+y)+h(2x-y)=4[h(x+y)+h(x-y)+6h(x)]-6h(y)\eqno(2.4)$$ for all $x,y \in X$. Thus it is enough to
prove that $h$ satisfies in $(2.4).$ Replacing $x$ and $y$ by $2x$ and
$2y$ in $(1.4),$ respectively, we obtain
$$f(2(2x+y))+f(2(2x-y))=4[f(2(x+y))+f(2(x-y))]+2[f(4x)-4f(2x)]-6f(2y)\eqno (2.5)$$
for all $x,y \in X$. Since $g(2x)=4g(x)$ for all $x \in X$,
$$f(4x)=20 f(2x)-64f(x)\eqno(2.6)$$ for all $x \in X$.
By $(2.5)$ and $(2.6),$ we get
$$f(2(2x+y))+f(2(2x-y))=4[f(2(x+y))+f(2(x-y))]+32[f(2x)-4f(x)]-6f(2y) \eqno(2.7)$$
for all $x,y \in X$. By  multiplying both sides of $(1.4)$ by $4$,
we get
$$4[f(2x+y)+f(2x-y)]=16[f(x+y)+f(x-y)]+8[f(2x)-4f(x)]-24f(y)$$
for all $x,y \in X$. If we subtract the last equation from
$(2.7),$ we obtain
\begin{align*}
h(2x+y)+h(2x-y)&=[f(2(2x+y))-4f(2x+y)]+[f(2(2x-y))-4f(2x-y)]\\
&=4[f(2(x+y))-4f(x+y)]+4[f(2(x-y))-4f(x-y)]\\
&+24[f(2x)-4f(x)]-6[f(2y)-4f(y)]\\
&=4[h(x+y)+h(x-y)+6h(x)]-6h(y) \hspace{4cm}
\end{align*}
for all $x,y \in X$.

Therefore, the mapping $h:X \rightarrow Y$ is
quartic. This completes the proof of the lemma.
\end{proof}

\begin{thm}\label{t2} A mapping $f:X \rightarrow Y$ satisfies
$(1.4)$ for all $x,y\in X$ if and only if there exist a unique
symmetric multi-additive mapping $M:X^4
\rightarrow Y$ and a unique symmetric bi-additive mapping $B:X
\times X \rightarrow Y$  such that
$$f(x)=M(x,x,x,x)+B(x,x)$$ for all $x\in X.$
\end{thm}
\begin{proof} Let $f$ satisfies $(1.4)$ and assume that
$g, h :X\rightarrow Y$ are mappings defined by
$$g(x):=f(2x)-16f(x), \hspace{2cm} h(x):=f(2x)-4f(x)\hspace{2cm}$$
 for all $x\in X.$ By Lemma 2.1, we obtain that the
 mappings $g$ and $h$ are quadratic and quartic, respectively, and
$$f(x)=\frac{1}{12}h(x)-\frac{1}{12}g(x)$$
for all $x\in X.$

Therefore, there exist a unique symmetric
multi-additive mapping $M:X^4 \rightarrow
Y$ and a unique symmetric bi-additive mapping $B:X\times
X\rightarrow Y$ such that $\frac{1}{12}h(x)=M(x,x,x,x)$ and
$\frac{-1}{12}g(x)=B(x,x)$ for all $x\in X$(see
\cite{Ac-Dh,Pa-Ba}). So
$$f(x)=M(x,x,x,x)+B(x,x)$$ for all $x\in X.$
 The proof of the converse is obvious.
\end{proof}
%================================================================
\section{Stability}
%\setcounter{equation}{0}
%================================================================
Throughout this section, assume that  $X$ is a real linear space and
$(Y,\mu,T)$ is a complete RN-space.
\begin{thm}\label{t2}
Let $f:X \to Y$ be a maping with $f(0)=0$ for which there is
$\rho:X \times X \to D^+$ ( $\rho(x,y)$ is denoted by $\rho_{x,y}$
) with the property:
$$\mu_{f(2x+y)+f(2x-y)-4f(x+y)-4f(x-y)-2f(2x)+8f(x)+6f(y)}(t) \geq \rho_{x,y}(t)\eqno(3.1)$$ for all $x,y \in X$ and all
$t>0.$ If
\begin{align*}
\lim_{n\to\infty}T_{i=1}^{\infty}(\rho_{2^{n+i-1}x,2^{n+i-1}x}(\frac{2^{2n+i}t}{4})&+\rho_{2^{n+i-1}x,
2.2^{n+i-1}x}(2^{2n+i}t)\\
&+\rho_{0,2^{n+i-1}x}(\frac{3.2^{2n+i}t}{4}))=1 \hspace{3.1cm}(3.2)
\end{align*}
and
$$\lim_{n \to \infty}\rho_{2^{n}x,2^{n}y}(2^{2n}t)=1 \eqno(3.3)$$
for all $x,y\in X$ and all $t>0$, then there exists a unique
quadratic mapping $Q_1:X \to Y$ such that
$$\mu_{f(2x)-16f(x)-Q_1(x)}(t) \geq T_{i=1}^{\infty}(\rho_{2^{i-1}x,2^{i-1}x}(\frac{2^{i}t}{4})+\rho_{2^{i-1}x,2.2^{i-1}x}(2^{i}t)+\rho_{0,2^{i-1}x}(\frac{3.2^{i}t}{4}))
\eqno(3.4)$$ for all $x \in X$ and all $t > 0.$
\end{thm}
\begin{proof}
Putting $y=x$ in $(3.1),$ we obtain
$$\mu_{f(3x)-6f(2x)+15f(x)}(t) \geq \rho_{x,x}(t) \eqno(3.5)$$
for all $x\in X.$ Letting $y=2x$ in $(3.1),$ we get
$$\mu_{f(4x)-4f(3x)+4f(2x)+8f(x)-4f(-x)}(t) \geq \rho_{x,2x}(t)\eqno(3.6)$$
for all $x \in X.$ Putting $x=0$ in $(3.1),$ we obtain
$$\mu_{3f(y)-3f(-y)}(t) \geq \rho_{0,y}(t) \eqno(3.7)$$
for all $y\in X.$ Replacing $y$ by $x$ in $(3.7),$ we see that
$$\mu_{3f(x)-3f(-x)}(t) \geq \rho_{0,x}(t) \eqno(3.8)$$
for all $x\in X.$ It follows from $(3.6)$ and $(3.8)$ that
$$\mu_{f(4x)-4f(3x)+4f(2x)+4f(x)}(t)\geq \rho_{x,2x}(t)+\rho_{0,x}(\frac{3t}{4})\eqno(3.9)$$ for all $x \in X.$
If we add $(3.5)$ to $(3.9),$ then we have
$$\mu_{f(4x)-20f(2x)+64f(x)}(t)\geq \rho_{x,x}(\frac{t}{4})+\rho_{x,2x}(t)+\rho_{0,x}(\frac{3t}{4}). \eqno(3.10)$$
Let
$$\psi_{x,x}(t)=\rho_{x,x}(\frac{t}{4})+\rho_{x,2x}(t)+\rho_{0,x}(\frac{3t}{4})\eqno(3.11)$$
for all $x \in X$. Then we get
$$\mu_{f(4x)-20f(2x)+64f(x)}(t)\geq \psi_{x,x}(t)\eqno(3.12)$$
for all $x \in X$ and all $t>0.$ Let $g:X \to Y$ be a mapping
defined by $g(x):=f(2x)-16f(x)$.  Then we conclude that
$$\mu_{g(2x)-4g(x)}(t)\geq \psi_{x,x}(t)\eqno(3.13)$$ for all $x \in X.$
Thus we have $$\mu_{\frac{g(2x)}{2^2}-g(x)}(t) \geq \psi_{x,x}(2^2
t)\eqno(3.14)$$ for all $x \in X$ and all $t>0.$ Hence
$$\mu_{\frac{g(2^{k+1}x)}{2^{2(k+1)}}-\frac{g(2^kx)}{2^{2k}}}(t)\geq \psi_{2^kx,2^kx}(2^{2(k+1)}t)\eqno(3.15)$$
for all $x \in X$ and all $k \in \Bbb N.$ This means that
$$\mu_{\frac{g(2^{k+1}x)}{2^{2(k+1)}}-\frac{g(2^kx)}{2^{2k}}}(\frac{t}{2^{k+1}})\geq \psi_{2^kx,2^kx}(2^{k+1}t)\eqno(3.16)$$
for all $x \in X,$ $t>0$ and all $k \in \Bbb N.$ By the triangle inequality, from $1
> \frac{1}{2}+\frac{1}{2^2}+\cdots+\frac{1}{2^n},$  it follows
\begin{align*}
\mu_{\frac{g(2^nx)}{2^{2n}}-g(x)}(t) \geq
T_{k=0}^{n-1}(\mu_{\frac{g(2^{k+1}x)}{2^{2(k+1)}}-\frac{g(2^kx)}{2^{2k}}}(\frac{t}{2^{k+1}}))&\geq
T_{k=0}^{n-1}(\psi_{2^kx,2^kx}(2^{k+1}t))\\
&=T_{i=1}^{n}(\psi_{2^{i-1}x,2^{i-1}x}(2^it))\hspace{1.26cm}(3.17)
\end{align*}
for all $x \in X$ and $t>0.$ In order to prove the convergence of
the sequence $\{\frac{{g(2^nx)}}{2^{2n}}\}$, we replace $x$ with
$2^{m}x$ in $(3.17)$ to obtain that
$$
\mu_{\frac{g(2^{n+m}x)}{2^{2(n+m)}}-\frac{g(2^mx)}{2^{2m}}}(t)
\geq T_{i=1}^{n}
(\psi_{2^{i+m-1}x,2^{i+m-1}x}(2^{i+2m}t)).\eqno(3.18)
$$
Since the right hand side of the inequality $(3.18)$ tends to $1$
as $m$ and $n$ tend to infinity, the sequence
$\{\frac{{g(2^nx)}}{2^{2n}}\}$ is a Cauchy sequence. Thus we
may define $Q_1(x)=\lim_{n \to \infty}\frac{{g(2^nx)}}{2^{2n}}$
for all $x \in X$. Now we show that $Q_1$ is a quadratic mapping.
Replacing $x,y$ with $2^nx$ and $2^ny$  in $(3.1), respectively,$ we get
$$\mu_{g(2x+y)+g(2x-y)-4g(x+y)-4g(x-y)-2g(2x)+8g(x)+6g(y)}(t) \geq \rho_{2^nx,2^ny}(2^{2n}t).\eqno(3.19)$$ Taking the limit as $n \to
\infty$, we find that $Q_1$ satisfies $(1.4)$ for all $x,y \in X$.
By Lemma 2.1, the mapping $Q_1:X \to Y$ is quadratic.

Letting the limit as $n \to \infty$ in $(3.17)$, we get $(3.4)$
by $(3.11).$

Finally, to prove the uniqueness of the quadratic
mapping $Q_1$ subject to $(3.4),$ let us assume that there exists
another quadratic mapping $Q_1'$ which satisfies $(3.4).$ Since
$Q_1(2^nx)=2^{2n}Q_1(x),$ $Q_1'(2^nx)=2^{2n}Q_1'(x)$ for all $x
\in X$ and $n \in \Bbb N,$ from $(3.4)$, it follows that
\begin{align*}
\mu_{Q_1(x)-Q_1'(x)}(2t)&=\mu_{Q_1(2^nx)-Q_1'(2^{n}x)}(2^{2n+1}t)\\
&\geq
T(\mu_{Q_1(2^nx)-g(2^nx)}(2^{2n}t),\mu_{g(2^nx)-Q_1'(2^nx)}(2^{2n}t))\\
&\geq
T(T_{i=1}^{\infty}(\rho_{2^{i+n-1}x,2^{i+n-1}x}(\frac{2^{2n+i}t}{4})+\rho_{2^{i+n-1}x,2.2^{i+n-1}x}(2^{2n+i}t)\\
&+\rho_{0,2^{i+n-1}x}(\frac{3.2^{2n+i}t}{4})),
T_{i=1}^{\infty}(\rho_{2^{i+n-1}x,2^{i+n-1}x}(\frac{2^{2n+i}t}{4})\\
&+\rho_{2^{i+n-1}x,2.2^{i+n-1}x}(2^{2n+i}t)+\rho_{0,2^{i+n-1}x}(\frac{3.2^{2n+i}t}{4})))
\hspace{2cm}(3.20)
\end{align*}
for all $x \in X$ and all $t > 0$. By letting $n \to \infty$ in
$(3.20),$ we conclude that $Q_1=Q_1'$.
\end{proof}
\begin{thm}\label{t2}
Let $f:X \to Y$ be a mapping with $f(0)=0$ for which there is
$\rho:X \times X \to D^+$ ( $\rho(x,y)$ is denoted by $\rho_{x,y}$
) with the property:
$$\mu_{f(2x+y)+f(2x-y)-4f(x+y)-4f(x-y)-2f(2x)+8f(x)+6f(y)}(t) \geq \rho_{x,y}(t)\eqno(3.21)$$ for all $x,y \in X$ and all
$t>0.$ If
\begin{align*}
\lim_{n\to\infty}T_{i=1}^{\infty}(\rho_{2^{n+i-1}x,2^{n+i-1}x}(\frac{2^{4n+3i}t}{4})&+\rho_{2^{n+i-1}x,2.2^{n+i-1}x}(2^{4n+3i}t)\\
&+\rho_{0,2^{n+i-1}x}(\frac{3.2^{4n+3i}t}{4}))=1\hspace{2.9cm}(3.22)
\end{align*}
and
$$\lim_{n \to \infty}\rho_{2^{n}x,2^{n}y}(2^{4n}t)=1 \eqno(3.23)$$
for all $x,y\in X$ and all $t>0$, then there exists a unique
quartic mapping $Q_2:X \to Y$ such that
$$\mu_{f(2x)-4f(x)-Q_2(x)}(t) \geq T_{i=1}^{\infty}(\rho_{2^{i-1}x,2^{i-1}x}(\frac{2^{3i}t}{4})+\rho_{2^{i-1}x,2.2^{i-1}x}(2^{3i}t)+\rho_{0,2^{i-1}x}(\frac{3.2^{3i}t}{4}))\eqno(3.24)$$ for all $x \in X$ and all $t > 0.$
\end{thm}

\begin{proof}
Putting $y=x$ in $(3.21),$ we obtain
$$\mu_{f(3x)-6f(2x)+15f(x)}(t) \geq \rho_{x,x}(t) \eqno(3.25)$$
for all $x\in X.$ Letting $y=2x$ in $(3.21),$ we get
$$\mu_{f(4x)-4f(3x)+4f(2x)+8f(x)-4f(-x)}(t) \geq \rho_{x,2x}(t)\eqno(3.26)$$
for all $x \in X.$ Putting $x=0$ in $(3.21),$ we obtain
$$\mu_{3f(y)-3f(-y)}(t) \geq \rho_{0,y}(t) \eqno(3.27)$$
for all $y\in X.$ Replacing $y$ by $x$ in $(3.27),$ we get
$$\mu_{3f(x)-3f(-x)}(t) \geq \rho_{0,x}(t) \eqno(3.28)$$
for all $x\in X.$ It follows from $(3.6)$ and $(3.28)$ that
$$\mu_{f(4x)-4f(3x)+4f(2x)+4f(x)}(t)\geq \rho_{x,2x}(t)+\rho_{0,x}(\frac{3t}{4})\eqno(3.29)$$ for all $x \in X.$
If we add $(3.25)$ to $(3.29),$ then we have
$$\mu_{f(4x)-20f(2x)+64f(x)}(t)\geq \rho_{x,x}(\frac{t}{4})+\rho_{x,2x}(t)+\rho_{0,x}(\frac{3t}{4}). \eqno(3.30)$$
Let
$$\psi_{x,x}(t)=\rho_{x,x}(\frac{t}{4})+\rho_{x,2x}(t)+\rho_{0,x}(\frac{3t}{4})\eqno(3.31)$$
for all $x \in X$. Then we get
$$\mu_{f(4x)-20f(2x)+64f(x)}(t)\geq \psi_{x,x}(t)\eqno(3.32)$$
for all $x \in X$ and all $t>0.$ Let $h:X \to Y$ be a mapping
defined by $h(x):=f(2x)-4f(x)$.  Then we conclude that
$$\mu_{h(2x)-16h(x)}(t)\geq \psi_{x,x}(t)\eqno(3.33)$$ for all $x \in X.$
Thus we have $$\mu_{\frac{h(2x)}{2^4}-h(x)}(t) \geq \psi_{x,x}(2^4
t)\eqno(3.34)$$ for all $x \in X$ and all $t>0.$ Hence
$$\mu_{\frac{h(2^{k+1}x)}{2^{4(k+1)}}-\frac{h(2^kx)}{2^{4k}}}(t)\geq \psi_{2^kx,2^kx}(2^{4(k+1)}t)\eqno(3.35)$$
for all $x \in X$ and all $k \in \Bbb N.$ This means that
$$\mu_{\frac{h(2^{k+1}x)}{2^{4(k+1)}}-\frac{h(2^kx)}{2^{4k}}}(\frac{t}{2^{k+1}})\geq \psi_{2^kx,2^kx}(2^{3(k+1)}t)\eqno(3.36)$$
for all $x \in X,$ $t>0$ and all $k \in \Bbb N.$ By the triangle inequality, from $1
> \frac{1}{2}+\frac{1}{2^2}+\cdots +\frac{1}{2^n},$  it follows
\begin{align*}
\mu_{\frac{h(2^nx)}{2^{4n}}-h(x)}(t) \geq
T_{k=0}^{n-1}(\mu_{\frac{h(2^{k+1}x)}{2^{4(k+1)}}-\frac{h(2^kx)}{2^{4k}}}(\frac{t}{2^{k+1}}))&\geq
T_{k=0}^{n-1}(\psi_{2^kx,2^kx}(2^{3(k+1)}t))\\
&=T_{i=1}^{n}(\psi_{2^{i-1}x,2^{i-1}x}(2^{3i}t))\hspace{1.25cm}(3.37)
\end{align*}
for all $x \in X$ and all $t>0.$ In order to prove the convergence of
the sequence $\{\frac{{h(2^nx)}}{2^{4n}}\}$, we replace $x$ with
$2^{m}x$ in $(3.37)$ to obtain that
$$
\mu_{\frac{h(2^{n+m}x)}{2^{4(n+m)}}-\frac{h(2^mx)}{2^{4m}}}(t)
\geq T_{i=1}^{n}
(\psi_{2^{i+m-1}x,2^{i+m-1}x}(2^{3i+4m}t)). \eqno(3.38)
$$
Since the right hand side of the inequality $(3.38)$ tends to $1$
as $m$ and $n$ tend to infinity, the sequence
$\{\frac{{h(2^nx)}}{2^{4n}}\}$ is a Cauchy sequence. Thus we
may define $Q_2(x)=\lim_{n \to \infty}\frac{{h(2^nx)}}{2^{4n}}$
for all $x \in X$. Now we show that $Q_2$ is a quartic mapping.
Replacing $x,y$ with $2^nx$ and $2^ny$ in $(3.21), respectively,$
we get
$$\mu_{h(2x+y)+h(2x-y)-4h(x+y)-4h(x-y)-2h(2x)+8h(x)+6h(y)}(t) \geq \rho_{2^nx,2^ny}(2^{4n}t).\eqno(3.39)$$ Taking the limit as $n \to
\infty$, we find that $Q_2$ satisfies $(1.4)$ for all $x,y \in X$.
By Lemma 2.1 we get that the mapping $Q_2:X \to Y$ is quartic.\\

Letting the limit as $n \to \infty$ in $(3.37)$, we get $(3.24)$
by $(3.31).$

 Finally, to prove the uniqueness of the quartic
mapping $Q_2$ subject to $(3.24),$ let us assume that there
exists a quartic mapping $Q_2'$ which satisfies $(3.24).$ Since
$Q_2(2^nx)=2^{4n}Q_2(x)$ and $Q_2'(2^nx)=2^{4n}Q_2'(x)$ for all $x
\in X$ and $n \in \Bbb N,$ from $(3.24)$, it follows that
\begin{align*}
\mu_{Q_2(x)-Q_2'(x)}(2t)&=\mu_{Q_2(2^nx)-Q_2'(2^{n}x)}(2^{4n+1}t)\\
&\geq T(\mu_{Q_2(2^nx)-h(2^nx)}(2^{4n}t),\mu_{h(2^nx)-Q_2'(2^nx)}(2^{4n}t))\\
&\geq T(T_{i=1}^{\infty}(\rho_{2^{i+n-1}x,2^{i+n-1}x}(\frac{2^{4n+3i}t}{4})+\rho_{2^{i+n-1}x,2.2^{i+n-1}x}(2^{4n+3i}t)\\
&+\rho_{0,2^{i+n-1}x}(\frac{3.2^{4n+3i}t}{4})),
T_{i=1}^{\infty}(\rho_{2^{i+n-1}x,2^{i+n-1}x}(\frac{2^{4n+3i}t}{4})\\
&+\rho_{2^{i+n-1}x,2.2^{i+n-1}x}(2^{4n+3i}t)+\rho_{0,2^{i+n-1}x}(\frac{3.2^{4n+3i}t}{4})))\hspace{1.25cm}(3.40)
\end{align*}
for all $x \in X$ and all $t > 0$. By letting $n \to \infty$ in
$(3.40),$ we get that $Q_2=Q_2'$.
\end{proof}

\begin{thm}\label{t2}
Let $f:X \to Y$ be a mapping with $f(0)=0$ for which there is
$\rho:X \times X \to D^+$ ( $\rho(x,y)$ is denoted by $\rho_{x,y}$
) with the property:
$$\mu_{f(2x+y)+f(2x-y)-4f(x+y)-4f(x-y)-2f(2x)+8f(x)+6f(y)}(t) \geq \rho_{x,y}(t)\eqno(3.41)$$ for all $x,y \in X$ and all
$t>0.$ If
\begin{align*}
&\lim_{n\to\infty}T_{i=1}^{\infty}(\rho_{2^{n+i-1}x,2^{n+i-1}x}(\frac{2^{4n+3i}t}{4})+\rho_{2^{n+i-1}
x,2.2^{n+i-1}x}(2^{4n+3i}t)+\rho_{0,2^{n+i-1}x}(\frac{3.2^{4n+3i}t}{4}))\\
&=1 \\ &
=\lim_{n\to\infty}T_{i=1}^{\infty}(\rho_{2^{n+i-1}x,2^{n+i-1}x}
(\frac{2^{2n+i}t}{4})+\rho_{2^{n+i-1}x,2.2^{n+i-1}x}(2^{2n+i}t)+\rho_{0,2^{n+i-1}x}
(\frac{3.2^{2n+i}t}{4}))\\
&\hspace{12.2cm}(3.42)
\end{align*}
and
$$\lim_{n \to \infty}\rho_{2^{n}x,2^{n}y}(2^{4n}t)=1=\lim_{n \to \infty}\rho_{2^{n}x,2^{n}y}(2^{2n}t) \eqno(3.43)$$
for all $x,y\in X$ and all $t>0$, then there exist a unique
quadratic mapping $Q_1:X \to Y$ and a unique quartic mapping
$Q_2:X \to Y$ such that
\begin{align*}
&\mu_{f(x)-Q_1(x)-Q_2(x)}(t)\\
&\geq T_{i=1}^{\infty}(\rho_{2^{i-1}x,2^{i-1}x}(3.2^{i}t)+\rho_{2^{i-1}x,2.2^{i-1}x}(12.2^{i}t)+\rho_{0,2^{i-1}x}(9.2^{i}t))\\
&+T_{i=1}^{\infty}(\rho_{2^{i-1}x,2^{i-1}x}(3.2^{3i})+\rho_{2^{i-1}x,2.2^{i-1}x}(12.2^{3i}t)+\rho_{0,2^{i-1}x}(9.2^{3i}))\hspace{2.5cm}(3.44)
\end{align*}
for all $x \in X$ and all $t > 0.$
\end{thm}
\begin{proof}
By Theorems 3.1 and 3.2, there exist a quadratic mapping $Q^{'}_1:X
\to Y$ and a  quartic mapping $Q^{'}_2:X \to Y$ such that
$$\mu_{f(2x)-16f(x)-Q^{'}_1(x)}(t) \geq T_{i=1}^{\infty}(\rho_{2^{i-1}x,2^{i-1}x}(\frac{2^{i}t}{4})+\rho_{2^{i-1}x,2.2^{i-1}x}(2^{i}t)+\rho_{0,2^{i-1}x}(\frac{3.2^{i}t}{4}))$$
and
$$\mu_{f(2x)-4f(x)-Q^{'}_2(x)}(t) \geq T_{i=1}^{\infty}(\rho_{2^{i-1}x,2^{i-1}x}(\frac{2^{3i}t}{4})+\rho_{2^{i-1}x,2.2^{i-1}x}(2^{3i}t)+\rho_{0,2^{i-1}x}(\frac{3.2^{3i}t}{4}))$$
for all $x \in X$ and all $t>0$. So it follows from the last
inequalities that
\begin{align*}
&\mu_{f(x)+\frac{1}{12}Q^{'}_1(x)-\frac{1}{12}Q^{'}_2(x)}(t)\\
&\geq T_{i=1}^{\infty}(\rho_{2^{i-1}x,2^{i-1}x}(3.2^{i}t)+\rho_{2^{i-1}x,2.2^{i-1}x}(12.2^{i}t)+\rho_{0,2^{i-1}x}(9.2^{i}t))\\
&+T_{i=1}^{\infty}(\rho_{2^{i-1}x,2^{i-1}x}(3.2^{3i})+\rho_{2^{i-1}x,2.2^{i-1}x}(12.2^{3i}t)+\rho_{0,2^{i-1}x}(9.2^{3i}))
\end{align*}
for all $x \in X$ and all $t>0$. Hence we obtain $(3.46)$ by letting
$Q_1(x)=-\frac{1}{12}Q^{'}_1(x)$ and
$Q_2(x)=\frac{1}{12}Q^{'}_2(x)$ for all $x \in X.$ The uniqueness
property of $Q_1$ and $Q_2,$ are trivial.
\end{proof}

\section{Acknowledgement}

The third author was supported by Korea Research Foundation Grant
funded by the Korean Government (KRF-2008-313-C00041).

{\small
%----------------------------------------------------------------------%
}
\end{document}